\documentclass{amsart}
\usepackage{verbatim}
\usepackage{amsmath}
\usepackage{amsthm}
\usepackage{amsfonts}
\usepackage{amssymb}

\def\0{{\bf0}}
\def\1{{\bf1}}
\def\2{{\bf2}}
\def\3{{\bf3}}

\def\SS{\tilde{\mathfrak{S}}}
\def\tF{\tilde{F}}
\def\W{\tilde{S}_n}
\def\tS{\tilde{S}}
\def\Z{\mathbb{Z}}
\def\a{\alpha}
\def\G{{\mathcal G}}
\def\B{{\mathcal B}}

\newtheorem{thm}{Theorem}
\newtheorem{cor}[thm]{Corollary}
\newtheorem{conj}[thm]{Conjecture}
\newtheorem{lem}[thm]{Lemma}
\newtheorem{prop}[thm]{Proposition}
\theoremstyle{definition}
\newtheorem{algorithm}[thm]{Algorithm}
\newtheorem{example}[thm]{Example}

\newtheorem{remark}[thm]{Remark}
\newtheorem{defi}[thm]{Definition}

\newcommand{\remind}[1]{$$\texttt{*** #1 ***}$$}
\date{January 2006}
\begin{document}
\author{Thomas Lam and Mark Shimozono}
\email{tfylam@math.harvard.edu} \email{mshimo@vt.edu}
\title{A Little bijection for affine Stanley symmetric functions}
\thanks{
M.S. was supported in part by NSF DMS-0401012.}

\begin{abstract}
Little~\cite{Lit} developed a combinatorial algorithm to study the
Schur-positivity of Stanley symmetric functions and the
Lascoux-Sch\"{u}tzenberger tree.  We generalize this algorithm to
affine Stanley symmetric functions, which were introduced recently
in~\cite{Lam}.
\end{abstract}
\maketitle

\section{Introduction}

A new family of symmetric functions, called {\it affine Stanley
symmetric functions} were recently introduced in~\cite{Lam}.  These
symmetric functions $\tF_w$, indexed by affine permutations $w \in
\tS_n$, are an affine analogue of the Stanley symmetric functions
$F_w$ which Stanley~\cite{Sta} introduced to enumerate the reduced
decompositions of a permutation $w \in S_n$.  Stanley symmetric
functions were later shown to be stable limits of the {\it Schubert
polynomials} ${\mathfrak S}_w$~\cite{LS,BJS}.  In the case that $w$
is a {\it Grassmannian permutation} the Stanley symmetric function
is equal to a Schur function. Shimozono conjectured and
Lam~\cite{Lam2} recently showed that the symmetric functions $\tF_w$
also had a geometric interpretation.  When $w$ is {\it affine
Grassmannian} then $\tF_w$ represents a Schubert class in the
cohomology $H^*({\mathcal G}/{\mathcal P})$ of the affine
Grassmannian and $\tF_w$ is called an {\it affine Schur function}.
Affine Schur functions were introduced by Lapointe and Morse
in~\cite{LM}, where they were called {\it dual $k$-Schur functions}.

A key property of the Stanley symmetric function $F_w$ is that its
expansion in terms of the Schur functions $s_\lambda$ involves
non-negative coefficients.  This was proved by Edelman and
Greene~\cite{EG} using an insertion algorithm and separately by
Lascoux and Sch\"{u}tzenberger~\cite{LS} via {\it transition
formulae} for Schubert polynomials.  These transition formulae lead
to a combinatorial object known as the {\it
Lascoux-Sch\"{u}tzenberger tree} which allows one to write a Stanley
symmetric function $F_w$ in terms of other Stanley symmetric
functions $F_v$ labeled by permutations ``closer to Grassmannian''.

Answering a question of Garsia, Little~\cite{Lit} recently gave a
combinatorial proof of the identity
\[
\sum_{\substack{u = v\cdot t_{r,s} \\ l(u)=l(v)+1}} F_u(X) =
\sum_{\substack{w = v \cdot t_{s',r} \\ l(w)=l(v)+1}} F_w(X),
\]
from which the Lascoux-Sch\"{u}tzenberger tree can be deduced.
Here $v \in S_n$ is any permutation, $t_{r,s}$ denotes a
transposition, and the summations are over $s > r$ and $s' < r$
respectively.

The aim of this article is to generalize Little's bijection and to
prove an affine analogue of the above identity for the affine
Stanley symmetric functions $\tF_w$.  Our techniques are mainly
combinatorial and the resulting bijection appears to be interesting
in itself -- for example, it seems to be closely related to the
affine Chevalley formula. Unfortunately, we have been unable to use
our affine Little bijection to prove that an arbitrary affine
Stanley symmetric function expands positively in terms of affine
Schur functions. This positivity follows from the results of
Lam~\cite{Lam2} combined with unpublished work of
Peterson~\cite{Pet}.


\section{Affine Stanley symmetric functions}
\subsection{Affine symmetric group}
Let $\W$ be the affine symmetric group. It is a Coxeter group with
simple reflections $\{s_i \mid i \in \Z/n\Z\}$ and relations
$s_i^2=1$ for all $i$, $(s_is_{i+1})^3=1$ for all $i$, and
$(s_is_j)^2=1$ for $i \neq j \pm 1$ not adjacent mod $n$.

One may realize $\W$ as the set of all bijections
$w:\Z\rightarrow\Z$ such that $w(i+n)=w(i)+n$ for all $i$ and
$\sum_{i=1}^n w(i) = \sum_{i=1}^n i$. In this realization, to
specify an element it suffices to give the ``window"
$[w(1),w(2),\dotsc,w(n)]$.

Given integers $r$ and $s$ such that $s\not\equiv r \mod n$, there
is a unique element $t_{r,s}\in \W$ such that, in the function
notation, $t_{r,s}(r)=s$ and $t_{r,s}(s)=r$, and $t_{r,s}(i)=i$ for
all $i$ such that $i \not\equiv r \mod n$ and $i\not\equiv s\mod n$.
We note that $t_{r,s}=t_{r',s'}$ if and only if there is an integer
$k$ such that $\{r,s\}=\{r'+kn,s'+kn\}$ as sets. In this notation
$s_i=t_{i,i+1}$ for all $i$.

For $a_i\in \Z/n\Z$ we call $a=a_1\dotsm a_l$ a reduced word for
$w\in \W$ and write $a\in R(w)$ if $w=s_{a_1}\dotsm s_{a_l}$ such
that $l$ is minimum.  We call $l = \ell(w)$ the length of $w$.
 Let $v \lessdot w$ denote the covering relation of the
\textit{strong Bruhat order} $\le$ on $\W$. By definition $v\lessdot
w$ if and only if there is a reflection $t_{r,s}\in \W$ such that
$w=vt_{r,s}$ and $\ell(w)=\ell(v)+1$.

Say that $w\in \W$ is a \textit{right $r$-cover} of $v\in \W$ if
$v\lessdot w$ with $w=v t_{r,s}$ for $r<s$. Say that $w$ is a
\textit{left $r$-cover} of $v$ if $v \lessdot w$ with $w= v t_{s,r}$
where $s<r$. Let $\Psi^+_r(v)$ denote the set of $r$-right covers of
$v$ and $\Psi^-_r(v)$ denote the set of $r$-left covers of $v$.

Now let $S_n \subset \W$ denote the symmetric group generated by
$s_1,s_2,\ldots,s_{n-1}$.  This is a parabolic subgroup of $\W$. The
minimal length coset representatives of $\W/S_n$ are called {\it
Grassmannian} and the set of such elements is denoted $\W^-$.

\subsection{Cyclically decreasing permutations}
Let $a = a_1 a_2 \cdots a_l$ be a reduced word.  Then $a$ is called
\emph{cyclically decreasing} if
\begin{enumerate}
\item
The multiset $A = \{a_1,a_2,\ldots,a_l\}$ is a set.
\item
If $i, i+1 \in A$ then $i+1$ occurs before $i$ in $a$, where indices
are considered modulo $n$. In particular, if $n-1,0\in A$ then $0$
appears before $n-1$ in $a$.
\end{enumerate}
A permutation $w\in\W$ is called \emph{cyclically decreasing} if
there is a reduced word for $w$ which is cyclically decreasing. Say
that a proper subset $I\subset \Z/n\Z$ is a \textit{cyclic interval}
if it has the form $I=\{i,i+1,i+2,\dotsc,i+j\}$ with indices taken
mod $n$. Let $(i+j)\dotsm (i+1)i$ be the word of $I$. The cyclically
decreasing elements of $\W$ are characterized as follows.

\begin{lem} \label{l:cycdecr}
\begin{enumerate}
\item Let $A\subset \Z/n\Z$ be a proper subset. Then a word with
underlying set $A$ is cyclically decreasing if and only if it is a
shuffle of the words of the maximal cyclic subintervals of $A$.
\item Every cyclically decreasing word with underlying set $A$ is a
reduced word for the same element $w(A)\in \W$.
\item Every reduced word for $w(A)$ is cyclically decreasing with
underlying set $A$.
\end{enumerate}
\end{lem}
\begin{proof} (1) holds by definition. Since reflections in
different maximal cyclic subintervals commute, it follows that every
cyclically decreasing word with underlying set $A$, is equivalent to
a single canonical word, namely, the concatenation of the words of
the maximal cyclic subintervals of $A$, with these words occurring
in decreasing order by first element. This given word is a reduced
word for some element $w(A)\in \W$, so all the cyclically decreasing
words with underlying set $A$ are as well. This proves (2). Since
there are no repeated reflections, the braid relations do not apply,
and the only equivalences among reduced words for $w(A)$ are
commutations between reflections in different maximal cyclic
subintervals of $A$. Since the set of reduced words is connected by
the Coxeter relations it follows that every reduced word for $w(A)$
is a shuffle of the prescribed sort. This proves (3).
\end{proof}

Lemma~\ref{l:cycdecr} has the following immediate consequence.

\begin{cor} \label{c:cycdecrBruhat} The strong Bruhat order on
the set of cyclically decreasing elements in $\W$ is isomorphic to
the boolean lattice on proper subsets of $\Z/n\Z$.
\end{cor}

Now let $w \in \tS_n$ of length $l = \ell(w)$ and suppose $\a =
(\a_1,\a_2,\ldots,\a_r)$ is a composition of $l$.  An
\emph{$\a$-decomposition} of $w$ is an ordered $r$-tuple of
cyclically decreasing affine permutations $(w^1,w^2,\ldots,w^r)\in
\tS_n^r$ satisfying $\ell(w^i) = \a_i$ and $w = w^1w^2\cdots w^r$.
The following definition is~\cite[Alternative Definition 2]{Lam}.

\begin{defi} The affine Stanley
symmetric function $\tF_w(X)$ is given by
\[
\tF_w(X) = \sum_{\a} (\mbox{number of $\a$-decompositions of
$w$})\cdot x^\a
\]
where the sum is over all compositions $\a$ of $\ell(w)$.
\end{defi}

It is shown in~\cite{Lam} that $\tF_w(X)$ is always a symmetric
function, though this fact will not be used in the current work.
When $w \in S_n \subset \tS_n$ is a normal permutation then, the
function $\tF_w$ is the usual Stanley symmetric function~\cite{Sta}.
When $w \in S_n \cap \W^-$ is a usual Grassmannian permutation, the
function $\tF_w$ is equal to some Schur function $s_\lambda$.

When $w \in \W^-$ is an affine Grassmannian permutation, then we say
that $\tF_w(X)$ is an {\it affine Schur function} or {\it dual
$k$-Schur function}.  There is a bijection $\theta: w
\leftrightarrow \lambda(w)$ between Grassmannian permutations $w \in
\W^-$ and partitions $\lambda(w)$ with no part greater than or equal
to $n$. The bijection $\theta$ sends a permutation $w$ with length
$l$ to a partition $\lambda$ with $l$ boxes.  We may thus label the
affine Schur functions by partitions $\tF_{\lambda(w)} := \tF_w$ so
that ${\rm deg}(\tF_{\lambda}) = |\lambda|$.  See~\cite{Lam} for
details.

\section{Affine Chevalley formula and affine Garsia-Little formula}
\subsection{Affine Chevalley formula}
Let $\G/\B$ denote the affine flag variety of type $A_{n-1}$;
see~\cite{Kum,KK}. The Bruhat decomposition of $\G$ induces a
decomposition of $\G/\B$ into Schubert cells \[ \G/\B = \cup_{w \in
\W} \Omega_w.\] Let $\SS_w \in H^*(\G/\B)$ denote the cohomology
class dual to $\Omega_w$.

The structure constants for the Schubert basis are denoted by
\begin{equation}
\label{eq:structure}
  \SS_u \SS_v = \sum_{w\in \W} c^w_{u,v} \SS_w
\end{equation}
 for $u,v\in\W$.  The following is a translation of the general Chevalley rule
(applicable to symmetrizable Kac-Moody groups) of Kostant and Kumar
\cite{KK} for the special case of $\W$.

\begin{prop} \cite{KK} For $v,w\in \W$ and any $r$, $c^w_{s_r,v}$ is
zero unless $w \gtrdot v$. In this case, writing $w = v t_{a,b}$
with $a<b$, $c^w_{s_r,v}$ is the number of times that $r$ occurs
modulo $n$ in the interval $[a,b-1]$.
\end{prop}

Conjecturally (\ref{eq:structure}) holds with $\SS_w$ replaced by
$\tF_w(X)$ everywhere.  Note that the functions $\tF_w(X)$ are not
linearly independent.  In particular we have the following
conjecture.

\begin{conj} \label{c:Chevalley}
\begin{equation} \label{e:tFChevalley}
  \tF_{s_r} \tF_v = \sum_{w\gtrdot v} c^w_{s_r,v} \tF_w.
\end{equation}
\end{conj}

This conjecture follows from an affine Schensted algorithm developed
in joint work with Lapointe and Morse~\cite{LLMS}.  It also follows
from unpublished geometric work of Peterson~\cite{Pet} and the
results in~\cite{Lam2}.

\begin{example} Let $n=4$. We use the window notation.
Let $v=[2,3,0,5]$ and $r=2$. The elements $w$ such that
$c_{s_r,v}^w$ nonzero are $w_1=[2,5,0,3]$ and $w_2=[2,4,-1,5]$. Now
$w_1 = v t_{2,4}$ and $2$ occurs once mod 4 in $\{2,3\}$, and
$w_2=vt_{2,7}$ and $2$ occurs twice mod 4 in $\{2,3,4,5,6\}$.
Therefore $\tF_{s_2} \tF_v = \tF_{w_1}+2\tF_{w_2}$.
\end{example}

\subsection{Affine Garsia-Little Formula}
The following, our main result, is an affine analogue of an identity
for Stanley functions observed by Garsia~\cite{Gar}, for which David
Little~\cite{Lit} found a combinatorial proof.

\begin{thm} \label{t:affLittle} For any $r\in\Z$ and $v\in \W$,
\begin{equation} \label{e:afflittle}
  \sum_{u\in \Psi_r^-(v)} \tF_u = \sum_{w\in \Psi_r^+(v)} \tF_w.
\end{equation}
\end{thm}

\begin{example}\label{e:positivity} Let $n=4$ and $v=[-1,1,4,6]=s_3s_1s_0$. For $r=1$ we
have $\Psi^+_1(v)=\{[1,-1,4,6]\}$ and
$\Psi^-_1(v)=\{[-3,3,4,6],[-2,1,4,7]\}$ and all the Chevalley
coefficients are 1:
\begin{equation*}
\tF_{[1,-1,4,6]}= \tF_{[-3,3,4,6]}+\tF_{[-2,1,4,7]}.
\end{equation*}
With the same $n$ and $v$ but $r=2$, we have
$\Psi^+_2(v)=\{[-1,4,1,6],[-3,3,4,6]\}$ and
$\Psi^-_2(v)=\{[1,-1,4,6],[-1,0,5,6]\}$ and all the Chevalley
coefficients are 1:
\begin{equation*}
\tF_{[-1,4,1,6]}+\tF_{[-3,3,4,6]}=
\tF_{[1,-1,4,6]}+\tF_{[-1,0,5,6]}.
\end{equation*}
Using these equations together one may find the equation
\begin{equation*}
  \tF_{[-1,4,1,6]} = \tF_{[-2,1,4,7]}+\tF_{[-1,0,5,6]}
\end{equation*}
The latter two are the affine Schur functions indexed by the
partitions $(2,1,1)$ and $(2,2)$ respectively (see~\cite{Lam}).
\end{example}

\begin{remark}
Theorem~\ref{t:affLittle} is a consequence of Conjecture
\ref{c:Chevalley}. Note that $\tF_{s_r}$ is the Schur function $s_1$
for all $r$. Subtracting \eqref{e:tFChevalley} for $\tF_{s_r} \tF_v$
and $\tF_{s_{r+1}} \tF_v$ one obtains \eqref{e:afflittle}.
\end{remark}

\subsection{Garsia-Little Formula and the Lascoux-Sch\"{u}tzenberger
Tree} Theorem~\ref{t:affLittle} holds, with a slight modification in
a special case (explained below), with all the affine objects
(affine permutations, affine transpositions, affine Stanley
symmetric functions) replaced by their usual $S_n$-counterparts. The
resulting {\it Garsia-Little formula} implies that a Stanley
symmetric function $F_w$ is Schur-positive, as follows.

Let $\sigma = \sigma_1 \sigma_2 \cdots \sigma_n \in S_n$ be a
permutation.  We set
\begin{eqnarray*}
r &=& \max(i \mid \sigma_i > \sigma_{i+1}), \\
s &=& \max(i
> r \mid \sigma_i < \sigma_r), \\
I &=& \{i < r \mid \sigma_i < \sigma_s \,\text{and for all $j \in
(i,r)$ we have $\sigma_j \notin (\sigma_i,\sigma_s)$}\}.
\end{eqnarray*}
Now let $\pi = \sigma \cdot t_{r,s}$.  One can check that we have
$\Psi^+_r(\pi) = \{ \sigma \}$ and $\Psi^-_r(\pi) = \{\pi \cdot
t_{i,r} \mid i \in I\}$.  If $I \neq \emptyset$ then
Theorem~\ref{t:affLittle} reads
\[
F_\sigma = \sum_{i \in I} F_{\pi \cdot t_{i,r}}.
\]
The permutations $\pi \cdot t_{i,r}$ are the children of $\sigma$ in
the Lascoux-Sch\"{u}tzenberger tree.  When $I = \emptyset$, the
corresponding equation fails to hold, but we declare $\sigma$ to
have a single child $1 \otimes \sigma = 1\, \sigma_1 \sigma_2 \cdots
\sigma_n \in S_{n+1}$ and we note that $F_\sigma = F_{1 \otimes
\sigma}$.  It can be shown that the process of (repeatedly) taking
children eventually results in permutations which are Grassmannian,
which are the leaves of the Lascoux-Sch\"{u}tzenberger tree. Since
the Stanley symmetric function indexed by $\sigma$ is equal to the
sum of those indexed by the children of $\sigma$ we conclude that it
is also the sum of those indexed by the leaves which are descendents
of $\sigma$.  A Stanley symmetric function indexed by a Grassmannian
permutation is a Schur function, so in particular every Stanley
symmetric function is Schur positive.

Unfortunately, a similar attempt to produce an ``affine
Lascoux-Sch\"{u}tzenberger tree'' fails because for an affine
permutation $w \in \tS_n$ there maybe no permutation $v$ and index
$r$ so that Theorem~7 involves only $\tF_w$ on one side.  However,
by solving simultaneous equations obtained from Theorem~7, we have
so far always been able to express an affine Stanley symmetric
function in terms of affine Schur functions, as demonstrated in
Example~\ref{e:positivity}.  It is our hope that the methods of this
paper will eventually lead to a combinatorial interpretation of the
coefficients $a_w^\lambda$ in the expansion $\tF_w = \sum_\lambda
a_w^\lambda \tF_\lambda$ of affine Stanley symmetric functions in
terms of affine Schur functions. These coefficients contain, for
example, the 3-point, genus zero, Gromov Witten invariants of the
Grassmannian, which are important numbers in combinatorics,
geometry, and representation theory; see~\cite{LM,Lam}.  The numbers
$a_w^\lambda$ also include as a special case the structure constants
for the multiplication of the homology of the affine Grassmannian in
the Schubert basis; see~\cite{Lam2,Pet}.

\section{Affine Little Bijection}

A \textit{$v$-marked word} is a pair $(a,i)$ where $a=a_1\dotsm a_l$
is a word with letters in $\{0,1,\dotsc,n-1\}$ and $i\in[1,l]$ is
the index of a distinguished letter in $a$ such that $a_1\dotsm
\hat{a_i}\dotsm a_l\in R(v)$.

We now define the \textit{affine Little graph}. It is a directed
graph whose vertices are the $v$-marked words. Given a $v$-marked
word $(a,i)$ with $a=a_1\dotsm a_l$, there is a unique directed edge
$(a,i)\rightarrow (a',j)$ leaving $(a,i)$, where $a'$ is the word
obtained from $a$ by replacing the letter $a_i$ by $a_i-1$ (mod $n$)
and the index $j$ is equal to $i$ if $a'$ is reduced and is
otherwise the unique index $j\not=i$ such that $a_1\dotsm \hat{a_j}
\dotsm a_l\in R(v)$, whose existence and uniqueness follows by Lemma
\ref{l:insert} (see Section~\ref{sec:coxeter}).

It is not hard to see that each vertex $(a,i)$ has a unique incoming
edge $(b,k)\rightarrow (a,i)$: if $a$ is reduced then $k=i$, and
otherwise, $k\not=i$ is the unique index such that $a_1\dotsm
\hat{a_k}\dotsm a_l\in R(v)$, and in either case, $b$ is obtained
from $a$ by replacing $a_k$ by $a_k+1$.

Since there are finitely many $v$-marked words, the connected
components of the affine Little graph are finite directed cycles. It
is clear from the definition that none of these cycles is a loop,
that is, it is never the case that $(a,i)\rightarrow (a,i)$.

Say that the $v$-marked word $(a,i)$ is \textit{reduced} if $a$ is a
reduced word for some $w\in \W$. Given a reduced $v$-marked word
$(a,i)$, write $\phi^v(a,i)=(b,j)$ where $(b,j)$ is the first
reduced $v$-marked word following $(a,i)$ on the cycle of the affine
Little graph containing $(a,i)$. The map $\phi^v$ defines a
bijection from the set of reduced $v$-marked words to itself. We
call an application of $\phi^v$ the \textit{affine Little
algorithm}.

Let $w\gtrdot v$, $w=vt_{r,s}$ and $a=a_1\dotsm a_l\in R(w)$. By
Proposition \ref{p:SEC} there exists a unique $i \in [1,l]$ so that
$a_1 \dotsm \hat{a_i} \dotsm a_l\in R(v)$ and $(a,i)$ is $v$-marked.
Therefore the set of reduced $v$-marked words is in bijection with
$\bigcup_{w\gtrdot v} R(w)$, and $\phi^v$ can be regarded as a
bijection from this set to itself.

\begin{example} Let $n=5$, $3410321042\in R(v)$, $r=2$,
and $b=3410\underline{2}321042\in \Psi^+_r(v)$ where the
distinguished reflection (the $j$-th) is underlined. The computation
of $\phi^v(b,j)$ is shown below. The indices $p(a,i)$ and $q(a,i)$
are those that appear in the proof of Theorem \ref{t:bijection}.
$\phi^v(b,j)$ is given by the last row. In Figure \ref{f:Little} the
edges of the affine Little graph for $v$ go from each $v$-marked
word to the one in the next row. We give some additional data used
in the proof of Theorem \ref{t:bijection} with $r=2$. We note that
for the last row, literally $p(a,i)=-1$ and $q(a,i)=7$. However
$-6<2$ and $t_{-1,7}=t_{-6,2}$. In other words, one should identify
the pairs $(p,q)$ and $(p',q')$ if there is a $k$ such that $
p'=kn+p$ and $q'=kn+q$.
\begin{figure}
\begin{equation*}
\begin{array}{|c|c|c|} \hline %
(a,i) & p(a,i) & q(a,i) \\ \hline %
3410\underline{2}321042 & 2 & 5 \\ \hline %
3410132104\underline{2} & 2 & 3 \\ \hline %
34\underline{1}01321041 & 2 & 1 \\ \hline %
340\underline{0}1321041 & 2 & 3 \\ \hline %
340\underline{4}1321041 & -6 & 2 \\ \hline %
\end{array}
\end{equation*}
\caption{\label{f:Little} Computation of $\phi^v$.}
\end{figure}
\end{example}

\begin{thm} \label{t:bijection}
The map $\phi^v$ restricts to a bijection
\[\phi^v_r: R(\Psi^+_r(v)) \longrightarrow R(\Psi^-_r(v)).\]
\end{thm}
\begin{proof} Due to the symmetry of left and right $r$-covers and
the bijectivity of $\phi^v$, it suffices to show that $\phi^v$ maps
$\Psi^+_r(v)$ into $\Psi^-_r(v)$.

Given a $v$-marked word $(a,i)$ with $a=a_1\dotsm a_l$, let $x$ and
$y$ be the elements of $\W$ with reduced words $a_1\dotsm a_{i-1}$
and $a_{i+1}\dotsm a_l$ respectively. Let $w=s_{a_1}\dotsm s_{a_l}$
and $t=a_i$. Then $w=xs_t y = xy (y^{-1}s_ty)=v t_{p,q}$ where
$y(p)=t$ and $y(q)=t+1$. Note also that if $s_t y > y$ (which occurs
if $a$ is reduced), then $p<q$, and if $s_t y < y$ then $p>q$. We
shall use the notation $p=p(a,i)$ and $q=q(a,i)$ to emphasize the
dependence on $(a,i)$.

Let $b\in R(\Psi^+_r(v))$ and $(b,j)$ the corresponding $v$-marked
word. Let $\phi^v(b,j)=(c,k)$ where $c\in R(u)$. It is enough to
show that $q(c,k)=r$, for if so, then since $c$ is reduced,
$p(c,k)<q(c,k)=r$ and $u\in \Psi^-_r(v)$ as desired.

To this end we show that for all vertices $(a,i)$ on the path in the
affine Little graph from $(b,j)$ to $(c,k)$ except $(c,k)$, that
$p(a,i)=r$.

Let $(a,i)$ be such a vertex and let $p=p(a,i)$ and $q=q(a,i)$.

For the base case $(a,i)=(b,j)$. Since $b\in R(w)$ we have $p<q$ and
since $w\in \Psi^+_r(v)$, $p=r$ as required.

For the induction step suppose $(a,i)$ satisfies $p=r$. Let
$(a,i)\rightarrow (a',i')$ in the affine Little graph. Write
$x',y',t',p',q'$ for the quantities associated with $(a',i')$.

By definition $a'$ is obtained from $a$ by replacing $a_i=t$ by
$a_i-1=t-1$. Let $w'=s_{a'_1}\dots s_{a'_l}$. Then $w'=x s_{t-1} y =
x y (y^{-1} s_{t-1} y) = v t_{p',q'}$. In particular, since $z s_m
z^{-1} = z t_{m,m+1} z^{-1}=t_{z(m),z(m+1)}$ for all $z\in \W$ and
$m$, we have $\{t-1,t\}=\{y(p'),y(q')\}$ as sets. Let $s$ be such
that $y(s)=t-1$. Since $t=y(p)=y(r)$ we have $\{r,s\}=\{p',q'\}$.

Suppose $a'$ is not reduced. By Lemma \ref{l:insert}, $i'\not=i$.
Suppose $i'<i$. Then $s_{t'}y'<y'$ and $s_{t-1}y>y$. It follows that
$q'<p'$ and $s<r$, so that $p'=r$ as desired. Suppose $i'>i$. Then
$s_{t'}y'>y'$ and $s_{t-1}y<y$, so that $p'<q'$ and $r<s$ and again
$p'=r$ as desired.

Otherwise let $a'$ be reduced. Then $(a',i')=(c,k)$, $i=i'$, and
$y'=y$. By the reduced-ness of $a'$, $s_{t-1}y>y$ and $s<r$. Again
by reduced-ness $p'<q'$. Therefore $q'=r$ as desired.
\end{proof}

\section{Generalized affine Little algorithm}
We now generalize the affine Little algorithm of the previous
section from reduced words to $\alpha$-decompositions.

\begin{lem}
\label{l:cyclic} Let $v\lessdot w$ both be cyclically decreasing.
Then there exists a cyclically decreasing $v$-marked reduced word
$a$ for $w$ which $\phi^v$ maps to a cyclically decreasing
$v$-marked reduced word, for the element $w'$ say. Furthermore the
element $w'$ is independent of the choice of $a$.
\end{lem}
\begin{proof} Let $w=w(A)$ and $v=w(A-\{i\})$. Let $j$ be maximal
such that $\{i,i-1,\dotsc,i-j\}\subset A$, with indices taken mod
$n$. Let $A'$ be the proper subset of $[0,n-1]$ obtained from $A$ by
replacing $i$ with $i-j-1$. Let $w'=w(A')$.

Now take any reduced word $a$ of $w$ and apply $\phi^v$ to the word
$w$ with the reflection $i$ marked. The application of $\phi^v$ to
$a$ replaces the subword $i (i-1)\dotsm (i-j)$ by
$(i-1)\dotsm(i-j-1)$, resulting in $a'$, say, with $i-j-1$ marked.
The only way that $a'$ is not cyclically decreasing is if $i-j-2\in
A$ and it appears to the left of $i-j$ in $a$. In this case $i-j-2$
and $i-j$ are in different maximal cyclically decreasing
subintervals of $A$. By Lemma \ref{l:cycdecr} there is an $a\in
R(w)$ with $i-j-2$ to the right of $i-j$. With such a choice of $a$,
by Lemma \ref{l:cycdecr} $a'$ is a reduced word for $w(A')$, which
depends only on $A'$, that is, only on $A$ and $i$.
\end{proof}

In particular, for $v\in \W$ cyclically decreasing, $\phi^v$ induces
a permutation of the set of cyclically decreasing covers of $v$.
Denote this map by $\phi^v$.


Let $v\in \W$ and $w\in \Psi^+_r(v)$. Let $\alpha$ be fixed. We will
describe an algorithm which takes as input an $\alpha$-decomposition
$w=w^1w^2\dotsm w^r$ of $w\in \Psi^+_r(v)$, and outputs an
$\alpha$-decomposition $x=x^1 \cdots x^r$ of an element
$x\in\Psi^-_r(v)$.

\begin{algorithm}
\label{a:extendedAffineLittle} Initialise $y:=w$ and $y^i:=w^i$. By
Proposition \ref{p:SEC} and Corollary \ref{c:cycdecrBruhat} we may
write $v = y^1 y^2 \cdots (y^i)' \cdots y^r$ where $(y^i)'$ is
obtained from $y^i$ by removing a single simple reflection. This
simple reflection is the same for any reduced word for $y^i$.
\begin{enumerate}
\item
Treating $y^i$ as a $(y^i)'$-marked word, apply $\phi^{(y^i)'}$ to
$y^i$.
\item
If $y^1y^2\cdots \phi^{(y^i)'}(y^i) \cdots y^r$ is an
$\alpha$-decomposition of some permutation (that is, it is
``reduced''), we terminate with this as the output. Otherwise by
Lemma \ref{l:insert} there is a unique index $j \neq i$ so that $$v
= y^1 y^2 \cdots (y^j)' \cdots \phi^{(y^i)'}(y^i) \cdots y^r$$ where
$(y^j)'$ is obtained from $y^j$ by removing a single simple
reflection. Replace the $\alpha$-decomposition $y = y^1y^2\cdots y^i
\cdots y^r$ by $y^1y^2\cdots \phi^{(y^i)'}(y^i) \cdots y^r$ and set
$i:= j$.
\item
Return to 1.
\end{enumerate}
\end{algorithm}

The fact that the algorithm is well-defined is clear from Lemma
\ref{l:cyclic}.

\begin{thm}
\label{t:extendedAffineLittle} Algorithm
\ref{a:extendedAffineLittle} is a bijection between
$\alpha$-decompositions of permutations in $\Psi^+_r(v)$ and
$\alpha$-decompositions of permutations in $\Psi^-_r(v)$.
\end{thm}
\begin{proof}
The fact that the output is an $\alpha$-decomposition of a
permutation in $\Psi^-_r(v)$ follows from the same argument as in
Theorem \ref{t:bijection}, and the fact that $\phi_r$ is a
bijection follows from the reversibility.
\end{proof}

\begin{proof}[Proof of Theorem~\ref{t:affLittle}.]
Since the coefficient $[x^\alpha]\tF_w$ of $x^\alpha$ in the affine
Stanley symmetric function $\tF_w$ is given by the number of
distinct $\alpha$-decompositions of $w$,
Theorem~\ref{t:extendedAffineLittle} proves
Theorem~\ref{t:affLittle}.
\end{proof}

\begin{remark}
Given an $\alpha$-decomposition of $w$, it is not clear whether
there is an initial choice of reduced word $a$ for $w$ so that the
affine Little algorithm applied to $a$ naturally gives the same
$\alpha$-decomposition as the generalized affine Little algorithm.
\end{remark}

\section{A Coxeter-theoretic result}
\label{sec:coxeter} Let $(W,S)$ be a Coxeter system, that is, $W$ is
a group generated by a set of \textit{simple reflections} $S$
subject only to relations of the form $s^2=1$ for all $s\in S$ and
for $s\not=s'\in S$, $(ss')^{m(s,s')}=1$ for $m(s,s')\ge2$. $W$ is
called a Coxeter group. A \textit{reflection} is by definition a
conjugate in $W$ of an element of $S$.

In this section, to avoid double subscripting, we shall write $s_i$
for an arbitrary simple reflection. A \textit{reduced word} for
$w\in W$ is a factorization $w=s_1\dotsm s_r$ of $w$ with $s_i\in S$
such that $r$ is minimum. The number $r$ is the length of $w$,
denoted $\ell(w)$.

Let $v \lessdot w$ denote the covering relation of the
\textit{strong Bruhat order} $\le$ on $W$. By definition $v\lessdot
w$ if and only if there is a reflection $t\in W$ such that $w=vt$
and $\ell(w)=\ell(v)+1$.

The goal of this section is to establish Lemma~\ref{l:insert}, which
is used crucially in the affine Little bijection.

\begin{prop} \label{p:SEC} (Strong Exchange Condition) \cite[Theorem 5.8]{Hum}
Let $w=s_1s_2\dotsm s_r$ with $s_i\in S$ not necessarily reduced.
Suppose $t$ is a reflection such that $\ell(wt)<\ell(w)$. Then there
is an index $i$ such that $wt = s_1\dotsm \hat{s_i} \dotsm s_r$
where the hatted reflection is omitted. Moreover if the expression
for $w$ is reduced then the index $i$ is unique.
\end{prop}

\begin{cor} \label{c:deletion} \cite[Cor. 5.8]{Hum} Suppose
$w=s_1\dotsm s_r$ with $s_i\in S$ and $r>\ell(w)$. Then there are
indices $i<j$ for which $w=s_1\dotsm \hat{s_i}\dotsm\hat{s_j}\dotsm
s_r$.
\end{cor}

\begin{lem} \label{l:strong} \cite[Lemma 5.11]{Hum} Let $w' \lessdot w$.
Suppose there is an $s\in S$ such that $w' < w's$ and $w's\not=w$.
Then both $w<ws$ and $w's<ws$.
\end{lem}

\begin{lem} \label{l:redhat} Let $w=s_1\dotsm s_r$ with $s_i\in S$
and $1\le i\le r$. Then there are unique reflections $t,t'$ such
that $wt=t'w=s_1\dotsm\hat{s_i}\dotsm s_r$.
\end{lem}
\begin{proof} Let $x=s_1\dotsm s_{i-1}$ and $y=s_{i+1}\dotsm s_r$.
Then $s_1\dotsm\hat{s_i}\dotsm s_r=xy$ and $w=xsy=xy(y^{-1}sy)$.
Letting $t=y^{-1}sy$ and $t'=xsx^{-1}$ we have $wt=xy$ and
$t'w=xsx^{-1}xsy=xy$ as desired.
\end{proof}

\begin{lem} \label{l:lengthadd} Let $x,y\in W$. Then
$\ell(xy)<\ell(x)+\ell(y)$ if and only if there is a reflection $t$
such that $xt<x$ and $ty<y$.
\end{lem}
\begin{proof} For the converse, if $t$ is a reflection such that
$xt<x$ and $ty<y$, then $\ell(xy)=\ell(xtty)\le \ell(xt)+\ell(ty)
\le \ell(x)-1+\ell(y)-1$ as desired.

Suppose $x,y\in W$ are such that $\ell(xy)<\ell(x)+\ell(y)$. Let
$x=s_1\dotsm s_l$ and $y=s'_1\dotsm s'_m$ be reduced. Then
$xy=s_1\dotsm s_l s'_1\dotsm s'_m$ is not reduced. By Corollary
\ref{c:deletion}, $xy$ has a factorization obtained by removing two
of the simple reflections. Suppose both of the reflections are
removed from the reduced word for $x$, that is, $xy=s_1\dotsm
\hat{s_i}\dotsm \hat{s_j}\dotsm s_l s'_1\dotsm s'_m$ for some $1\le
i<j\le l$. Multiplying both sides by $y^{-1}=s'_m\dotsm s'_1$ we
have $x=s_1\dotsm \hat{s_i}\dotsm \hat{s_j}\dotsm s_l$ which
contradicts the assumption that $x=s_1\dotsm s_l$ was reduced.
Similarly both reflections cannot be removed from the reduced word
for $y$. Therefore there are indices $1\le i\le l$ and $1\le j\le m$
such that $xy=s_1\dotsm\hat{s_i}\dotsm s_l s'_1\dotsm
\hat{s'_j}\dotsm s'_m$. By Lemma \ref{l:redhat} there are
reflections $t,t'\in W$ such that $xt=s_1\dotsm\hat{s_i}\dotsm s_l$
and $t'y=s'_1\dotsm \hat{s'_j}\dotsm s'_m$. Therefore
$xtt'y=s_1\dotsm \hat{s_i}\dotsm s_l s'_1\dotsm \hat{s'_j}\dotsm
s'_m=xy$. It follows that $tt'=1$ and $t=t'$ since $t,t'$ are
reflections. Since $xt$ admits a shorter factorization in simple
reflections that $x$ does, $xt<x$. Similarly $ty=t'y<y$ as desired.
\end{proof}

\begin{lem} \label{l:insert} Suppose $s_i\in S$ for $1\le i\le r$
is such that $v=s_1\dotsm \hat{s_i}\dotsm s_r$ is reduced, where
$\hat{s_i}$ means that the factor $s_i$ is removed. Suppose also
that $w=s_1\dotsm s_r$ is not reduced. Then there is a unique
index $j\not=i$ such that $s_1\dotsm \hat{s_j}\dotsm s_r$ is
reduced. Moreover $v=s_1\dotsm \hat{s_j}\dotsm s_r$.
\end{lem}
\begin{proof} Let $x=s_1\dotsm s_{i-1}$ and $y=s_{i+1}\dotsm s_r$.
Since $v=xy=s_1\dotsm \hat{s_i}\dotsm s_r$ is reduced, both of the
expressions for $x$ and $y$ are reduced and
$\ell(xy)=\ell(x)+\ell(y)$. By Lemma \ref{l:lengthadd} for $x,y\in
W$ and the reflection $s_i$, either $xs_i>x$ or $s_iy>y$. Without
loss of generality we assume that $xs_i>x$; one may reduce to this
case by reversing the reduced words and replacing elements of $W$ by
their inverses.

By assumption $xs_i=s_1\dotsm s_i$ is reduced. If $s_1\dotsm s_r$ is
reduced then we are done. Suppose not. Let $j$ be the minimum index
such that $i+1\le j\le l$ and $s_1\dotsm s_j$ is not reduced. Then
$s_1\dotsm s_{j-1} \gtrdot s_1\dotsm \hat{s_i} \dotsm s_{j-1}$.
Applying Lemma \ref{l:strong} to this covering relation and the
simple reflection $s_j$, we have either (a)
$s_1\dotsm\hat{s_i}\dotsm s_j=s_1\dotsm s_{j-1}$ or (b) $s_1\dotsm
s_j \gtrdot s_1\dotsm \hat{s_i}\dotsm s_j$. If (b) holds, since
$s_1\dotsm\hat{s_i}\dotsm s_j$ is reduced it follows that $s_1\dotsm
s_j$ is also, contradicting our choice of $j$. So (a) must hold.
Right multiplying by $s_{j+1}\dotsm s_r$ we have
$v=s_1\dotsm\hat{s_i}\dotsm s_r = s_1\dotsm \hat{s_j}\dotsm s_r$.
But $s_1\dotsm \hat{s_i}\dotsm s_r$ is reduced and $s_1\dotsm
\hat{s_j}\dotsm s_l$ has the same number of reflections, so it too
must be reduced.

For uniqueness, suppose there is a $k$ distinct from $i$ and $j$
such that $v=s_1\dotsm \hat{s_k}\dotsm s_r$ is reduced. We now treat
the three indices $i,j,k$ interchangeably and suppose without loss
of generality that $i<j<k$. Using the reduced words for $v$ that
omit $s_i$ and $s_j$, we have $\hat{s_i}\dotsm s_j = s_i\dotsm
\hat{s_j}$. Right multiplying by $s_j$ we have $\hat{s_i}\dotsm
\hat{s_j} = s_i\dotsm s_j$. But $s_i\dotsm s_j$ is reduced, being a
subword of the reduced word $s_1\dotsm \hat{s_k}\dotsm s_r$. This is
a contradiction.
\end{proof}

\end{document}